\documentclass[12pt]{iopart}
\usepackage{amssymb}
\usepackage{amsfonts}
\usepackage{epsfig}
\usepackage{graphics}
\usepackage{iopams}

\newtheorem {theorem}{Theorem}

\newtheorem {definition}{Definition}

\newtheorem {lemma}{Lemma}

\newtheorem {remark}{Remark}

\newenvironment {proof}[1][Proof]{\noindent \textbf {#1.} }{\ \rule {0.5em}{0.5em}}
\begin{document}

\title{Geometric shadowing in slow-fast Hamiltonian systems}

\author{Niklas Br\"{a}nnstr\"{o}m$^{1}$ \footnote{The research of NB was funded by the Academy of Finland and the EU research
training network CODY.} , Emiliano de Simone$^{1}$ \footnote{
Funded by the Academy of Finland.} , Vassili Gelfreich$^{2}$}
\address{$^{1}$Department of Mathematics,
P.O. Box 68, FIN-00014 University of Helsinki,
Finland}
\address{$^{2}$Mathematics Institute,
University of Warwick,
CV4 7AL Coventry,
United Kingdom }

\date{}

\begin{abstract}
We study a class of slow-fast Hamiltonian systems with any finite number of
degrees of freedom, but with at least one slow one and two fast ones. At $%
\varepsilon =0$ the slow dynamics is frozen. We assume that the frozen
system (i.e. the unperturbed fast dynamics) has families of hyperbolic
periodic orbits with transversal heteroclinics.

For each periodic orbit we define an action $J.$ This action may be viewed
as an action Hamiltonian (in the slow variables). It has been shown in \cite%
{BG2007} that there are orbits of the full dynamics which shadow any \emph{%
finite }combination of forward orbits of $J$ for a time $t=O(\varepsilon
^{-1})$.

We introduce an assumption on the mutual relationship between the actions $%
J. $ This assumption enables us to shadow any continuous curve (of arbitrary
length) in the slow phase space for any time. The slow dynamics shadows the
curve as a purely geometrical object, thus the time on the slow dynamics has
to be reparameterised.
\end{abstract}

\maketitle

\section{Introduction}

Let us consider the Hamiltonian system defined by the Hamiltonian function 
\begin{equation*}
H=H(x,y,u,v;\varepsilon )
\end{equation*}%
and the symplectic form 
\begin{equation*}
\Omega =dy\wedge dx+\frac{1}{\varepsilon }dv\wedge du
\end{equation*}%
where $x,y\in {{\mathbb{R}}}^{2m}$ and $u,v\in {{\mathbb{R}}}^{2d}.$ We
assume that the parameter $\varepsilon $ is small, hence the Hamiltonian
system is "slow-fast" with $\left( x,y\right) $ being fast variables and $%
\left( u,v\right) $ being slow ones. This can readily be seen from the
equations of motion 
\begin{equation}
\left\{ 
\begin{array}{c}
\dot{x}=\frac{\partial H}{\partial y}, \\ 
\\ 
\,\dot{y}=-\frac{\partial H}{\partial x},%
\end{array}%
\right. \ {{\textrm{\ \ }}}%
\begin{array}{c}
\dot{u}=\varepsilon \frac{\partial H}{\partial v}, \\ 
\\ 
\dot{v}=-\varepsilon \frac{\partial H}{\partial u}.%
\end{array}
\label{eq:FullDyn}
\end{equation}%
We note that the form of $\left( \ref{eq:FullDyn}\right) $ is not unusual,
it appears in many applications.

If we set $\varepsilon =0$ the equations of motion are reduced to 
\begin{equation}
\left\{ 
\begin{array}{c}
\dot{x}=\frac{\partial H}{\partial y}, \\ 
\\ 
\,\dot{y}=-\frac{\partial H}{\partial x},%
\end{array}%
\right. \ {{\textrm{\ \ }}}%
\begin{array}{c}
\dot{u}=0,\bigskip \\ 
\\ 
\dot{v}=0.%
\end{array}
\label{eq:Frozen}
\end{equation}%
We refer to (\ref{eq:Frozen}) as the \emph{frozen system}.

We are interested in describing the slow dynamics of system (\ref{eq:FullDyn}%
). In general it is a difficult task to exactly describe the slow dynamics
and therefore often only approximations, or averaged solutions, are sought.
The procedure to obtain the averaged solutions, i.e. an averaging method, is
typically tailored for the particular class of systems it is applied to. The
averaging method for the classical problem where there is only one fast
degree of freedom is described in \cite{BM1961}. This method has been
generalised to systems with several fast degrees of freedom where the fast
system rotates with a constant vector of frequencies, see \cite{Neish1976}.
There are also averaging methods for slow-fast systems where the fast system
is uniformly hyperbolic, see \cite{Anosov1960} (or \cite{LM1988} for a
description in English). Contrary to \cite{Neish1976} the class of problems
studied in \cite{Anosov1960} contains systems which are fully coupled, i.e.
the fast system depends on the slow variables as well as the fast ones. Or,
more generally, if the fast system is ergodic and satisfies assumptions on
how fast time averages converge to space averages then, see \cite{Kifer2004}%
, the fast dynamics can be averaged out. Regardless of which averaging
method is employed the aim is to derive effective equations for the slow
dynamics which are independent of the fast variables $\left( x,y\right) .$
We may write such averaged slow dynamics as 
\begin{equation}
\dot{u}=\varepsilon \left\langle \frac{\partial H}{\partial v}\right\rangle
,\quad \dot{v}=-\varepsilon \left\langle \frac{\partial H}{\partial u}%
\right\rangle ,  \label{eq:Averaged}
\end{equation}%
where $\left\langle \cdot \right\rangle $ denotes averaging with respect to $%
(x,y)$. A crucial point in the justification of all averaging methods is to
verify that time averages can be approximated by space averages. This
property holds for most but not all trajectories. For example, if the fast
dynamics spends most of its time close to a periodic orbit there is, in
general, no reason to believe that the time average would converge to the
space average.

In fact when the fast dynamics spend most of its time close to periodic
orbits the slow dynamics of (\ref{eq:FullDyn}) may behave very differently
from the averaged dynamics (\ref{eq:Averaged}). In \cite{BG2007} it was
shown that if the system satisfies assumptions [\textbf{A1}] and [\textbf{A2}%
] (defined later) then there exist trajectories of the full dynamics whose
slow component deviates significantly from (\ref{eq:Averaged}). In
particular it was shown that there are orbits shadowing \emph{accessible
paths} of \emph{finite }length composed of forward trajectories of the
auxiliary systems 
\begin{equation*}
\dot{u}=\varepsilon \frac{\partial J_{c}}{\partial v},\quad \dot{v}%
=-\varepsilon \frac{\partial J_{c}}{\partial u},
\end{equation*}%
where $J_{c}$ is an action defined on a periodic orbit (labelled by $c$) in
the fast phase space (see Definition \ref{Def:Action}). This shadowing
result is valid on time-scales of order $O(\varepsilon ^{-1}).$ We note that
this method yields trajectories which deviate at the rate $O(\varepsilon )$
from (\ref{eq:Averaged}).

This approach is a generalisation of the mechanism proposed in \cite{GT2007}
for studying drift of the energy in a Hamiltonian system which depends on
time explicitly and slowly. In this set-up it was shown that switching
between fast periodic orbits does indeed provide the fastest possible rate
of energy growth in several situations (see \cite{GT2007}). An interesting
direct application of this theory is a rigorous proof of Fermi acceleration
for a class of billiards with slowly moving boundary, \cite{GT2008}.

The main result of this paper eliminates the upper bound on the time for
which the shadowing result holds, that is we provide a description of the
shadowing orbits for all times. Moreover, we show that for any continuous
curve in the slow phase space there is a trajectory of the full dynamics
whose slow component shadows it. We achieve this by refining the mechanism
in \cite{GT2007} and \cite{BG2007} such that for any two $O(\varepsilon )$
close points in the slow phase space there is a trajectory which starts in a
neighbourhood of the first point and ends up in a neighbourhood of the
second one. To prove this we have to introduce an assumption on the mutual
relationship between the actions $J_{c}$. This assumption ensures that any
vector in the tangent space of a point in the slow phase space can be
written as a linear combination of gradients of the actions $J_{c}$ where
all coefficients are positive. We identify this linear combination with a
path in the slow phase space and denote it \emph{guiding path }(which is a
generalisation of accessible path). We only consider guiding paths of length 
$O(\varepsilon )$ and show that for any such guiding path there exists a
trajectory of the full dynamics whose slow component shadows it. Thus by
approximating any given curve by a set of points which are $O(\varepsilon )$
apart we can find trajectories of the full dynamics whose slow component
shadows a guiding path between these points. Then \cite{GT2007} implies that
that there exists a trajectory which lies close to the union of trajectories
shadowing the guiding paths, hence there is a trajectory which shadows the
entire curve. The idea of the proof is similar to that of \cite{BG2007},
that is, using that the full dynamics takes place on normally hyperbolic
manifolds where certain action Hamiltonians are preserved for long times.

The paper is structured as follows. In the next section we present the
assumptions on the frozen system (\ref{eq:Frozen}) and state our main
result. Before proving this result we present an example where the
assumptions on the frozen system are verified to hold. In Section \ref%
{Sec:NormHyp} we provide a short summary of the results from \cite{BG2007}
and \cite{GT2007} on normal hyperbolicity which we will require to prove our
main theorem. In Section \ref{sec:ProofMainThm} we state two lemmas which we
then combine to prove Theorem \ref{Theorem:NewMain}. The proofs of the two
lemmas are of more technical nature and have been postponed to Section \ref%
{Sec:Technical} to increase the readability of the previous section.

\section{Set-up and Statement of the Result}

We impose a number of assumptions on the frozen system (\ref{eq:Frozen}).
Let $D\subset {{\mathbb{R}}}^{2d}$ be an open and bounded subset.

\begin{enumerate}
\item[{[\textbf{A1}]}] We assume that the frozen system (\ref{eq:Frozen}) has 
$n$ families of hyperbolic periodic orbits $L_{c}(u,v),$ $c\in \left\{
c_{1},\ldots ,c_{n}\right\} ,$ defined for all $\left( u,v\right) \in D.$

\item[{[\textbf{A2}]}] We assume that each of the periodic orbits has a
family of heteroclinic orbits to every other periodic orbit, i.e. for all $%
c_{i},c_{j}\in \left\{ c_{1},\ldots ,c_{n}\right\} $ and $\left( u,v\right)
\in D$ there are a pair of transversal heteroclinic orbits 
\begin{eqnarray*}
\Gamma _{c_{i}c_{j}}(u,v) &\subset &W^{u}(L_{c_{i}}(u,v))\cap
W^{s}(L_{c_{j}}(u,v)), \\
\Gamma _{c_{j}c_{i}}(u,v) &\subset &W^{u}(L_{c_{j}}(u,v))\cap
W^{s}(L_{c_{i}}(u,v)).
\end{eqnarray*}
\end{enumerate}

We note that under assumptions [\textbf{A1}] and [\textbf{A2}] the frozen
system has a family of uniformly hyperbolic invariant transitive sets $%
\Lambda _{(u,v)},$ also known as Smale horseshoes (see for instance \cite%
{Moser1973}). The dynamics on the Smale horseshoe can be described by
symbolic dynamics. Let $\Lambda :=\cup _{(u,v)\in D}\Lambda _{(u,v)}.$ For
each family of periodic orbits we define a family of actions.

\begin{definition}
\label{Def:Action}The action $J_{c}$ of a periodic orbit $L_{c}$ is defined
by the integral 
\begin{equation*}
J_{c}(u,v):=\oint_{L_{c}(u,v)}ydx.
\end{equation*}
\end{definition}

The function $J_{c}(u,v)$ is independent of the fast variables and can be
considered as a Hamiltonian function which generates some dynamics in the
slow variables 
\begin{equation}
\dot{u}=\frac{1}{T_{c}(u,v)}\frac{\partial J_{c}}{\partial v},\quad \dot{v}=-%
\frac{1}{T_{c}(u,v)}\frac{\partial J_{c}}{\partial u},  \label{eq:Guiding}
\end{equation}%
where $T_{c}$ is the period of the periodic orbit $L_{c}.$ Following the
notation of \cite{SVM2007} we refer to (\ref{eq:Guiding}) as the \textit{%
guiding system} and its equations of motion can be written in the concise
form 
\begin{equation}
\dot{z}=(\dot{u},\dot{v})=:X_{c}(z),  \label{eq:Concise}
\end{equation}%
where $X_{c}$ is the \emph{guiding vector field} .

\begin{enumerate}
\item[{[\textbf{A3}]}] We assume that there is a closed subset $\overline{{{%
\mathcal{D}}}}\subset D\subset {{\mathbb{R}}}^{2d}$ such that for all $%
(u,v)\in \overline{{{\mathcal{D}}}},$ \textbf{0} is inside the convex
envelope (or convex hull) of $\left\{ \nabla _{(u,v)}J_{c_{i}}\right\} $
where $i=1,\ldots ,n$.
\end{enumerate}

\begin{remark}
Although $n=2d+1$ actions are sufficient to satisfy \textbf{[A3]} we allow
for a larger set of actions to be considered since it may increase the size
of the domain $\overline{{{\mathcal{D}}}}.$ Nevertheless, at any point in $%
\overline{{{\mathcal{D}}}}$ we will use only $2d+1$ vectors (but not
necessarily the same $2d+1$ vectors for every point in $\overline{{{\mathcal{%
D}}}}$).

\end{remark}%

We can now state the main result. Let $\pi :\mathbb{R}^{2m+2d}\rightarrow 
\mathbb{R}^{2d}$ be the projection on the slow variables.

\begin{theorem}
\label{Theorem:NewMain}Let $\gamma :{\mathbb{R}}^{+}\rightarrow \overline{{%
\mathcal{D}}}$ be any continuous curve$.$ Assume that the frozen system (\ref%
{eq:Frozen}) satisfies assumptions \textbf{[A1],[A2]} and \textbf{[A3]}$.$
Then there exist positive constants $C$ and $\varepsilon _{0}$ such that for
any $0<\varepsilon <\varepsilon _{0}$ there exists a trajectory $\psi (t)$
of the full system (\ref{eq:FullDyn}) and a continuous monotone
reparameterisation of time ${\mathcal{T}}(t)$ such that the slow component $%
z(t):=\pi \psi (t)$ of $\psi (t)$ satisfies $z(0)=\gamma (0)$ and 
\begin{equation*}
\left\Vert z({\mathcal{T}}(t))-\gamma (t)\right\Vert \leq C\varepsilon ,
\end{equation*}%
for all $t\geq 0.$
\end{theorem}

Before proving the theorem we construct an example where assumptions \textbf{%
[A1], [A2]} and \textbf{[A3]} are shown to be satisfied.

\subsection{Example}
Let $\varepsilon =0$ and let $\mu $ be a small parameter. Consider a
Hamiltonian function of the form%
\[
H(x,y)=H_{0}(x,y)-\mu H_{1}(x,y)
\]%
where $(x,y)\in \mathbb{R}^{4}.$ Assume that the dynamics of the
"unperturbed" Hamiltonian $H_{0}$ has three hyperbolic periodic orbits $%
L_{c},$ $c\in \left\{ c_{1},c_{2},c_{3}\right\} $ with transversal
heteroclinic connections (hence assumption [\textbf{A1}] and [\textbf{A2}]
are satisfied). The implicit function theorem implies that the hyperbolic
periodic orbits persists for $\mu $ sufficiently small, and furthermore the
periodic orbits depend smoothly on $\mu .$ Therefore we can expand the
periodic orbit $L_{c}$ as well as the action $J_{c}$ in a power series in $%
\mu $. For the action we write%
\[
J_{c}=J_{0}^{c}+\mu J_{1}^{c}+O(\mu ^{2}),
\]%
where $J_{0}$ is the action of the periodic orbits of the "unperturbed"
Hamiltonian $H_{0}.$ $J_{0}$ and $J_{1}$ are constants as the Hamiltonian is
independent of $(u,v)$. \ The action $J_{0}$ and its first order correction $%
J_{1}$ are generically non-zero. Now, by abusing the notation slightly we
let $\mu $ depend on $(u,v)$ in the following way%
\begin{eqnarray*}
\mu =\mu (x,y,u,v)= &&\mu \sum_{i=1}^{3}\chi _{L_{c_{i}}+\delta
}(x,y)\varphi ^{c_{i}}(u,v) \\
&&+\mu \sum_{i=1}^{3}\chi \left( _{L_{c_{i}}+\delta }\right) ^{C}(x,y)f(x,y),
\end{eqnarray*}%
where $\mu $ is a small parameter, $\chi $ is the indicator function, $\
L_{c_{i}}+\delta $ is a $\delta $ neighbourhood of $L_{c_{i}},$ i.e.%
\[
L_{c_{i}}+\delta :=\left\{ (x,y)\textrm{ s.t. }dist((x,y),L_{c_{i}})\leq
\delta \right\} 
\]%
and $f$ is chosen as to interpolate $H_{1}$ into a $\mathcal{C}^{\infty }$
function. Then the action of the periodic orbit $L_{c}$ becomes%
\[
J_{c}(u,v)=J_{0}^{c}+\mu \varphi ^{c_{i}}(u,v)J_{1}^{c_{i}}+O(\mu ^{2}).
\]%
For example, by choosing%
\begin{eqnarray*}
\varphi ^{c_{1}} &=&sgn(J_{1}^{c_{1}})v, \\
\varphi ^{c_{2}} &=&-sgn(J_{1}^{c_{2}})u, \\
\varphi ^{c_{3}} &=&sgn(J_{1}^{c_{3}})u-sgn(J_{1}^{c_{3}})v,
\end{eqnarray*}%
where $sgn(A):=1$ if $A>0$ and $sgn(A):=-1$ if $A<0,$ will generate%
\begin{eqnarray*}
\nabla _{(u,v)}J^{c_{1}} &=&\mu \overrightarrow{e}_{u}, \\
\nabla _{(u,v)}J^{c_{2}} &=&\mu \overrightarrow{e}_{v}, \\
\nabla _{(u,v)}J^{c_{3}} &=&-\mu \overrightarrow{e}_{u}-\mu \overrightarrow{e%
}_{v},
\end{eqnarray*}%
This choice ensures that assumption [\textbf{A3] }is satisfied.

\section{Normal Hyperbolicity and symbolic dynamics\label{Sec:NormHyp}}

It was shown in \cite{BG2007} and \cite{GT2007} that the dynamics of the
full system $\left( \ref{eq:FullDyn}\right) $ can be studied using symbolic
dynamics. Here we provide a short summary of the results therein which we
will need to prove Theorem \ref{Theorem:NewMain}. For each periodic orbit $%
L_{c}$ we denote by $\Sigma _{c}$ a Poincar\'{e} section, and by $%
x_{i},y_{i},z_{i}$ (where $z_{i}=(u_{i},v_{i})$) we denote an intersection
of a solution of (\ref{eq:FullDyn}) with one of these Poincar\'{e} sections.
Let $\xi =\left\{ \xi _{i}\right\} _{i=-\infty }^{i=\infty }$ be a
bi-infinite sequence of letters $\xi _{i}\in \left\{ c_{1},\ldots
,c_{n}\right\} $ where $c_{i}$ is the index of the $n$ actions. The sequence 
$\xi $ is called the code of the dynamics. It has been shown that by
specifying an initial condition for the dynamics the code $\xi $ generates a
trajectory of (\ref{eq:FullDyn}) whose dynamics on the Poincar\'{e} sections
are given by 
\begin{equation}
z_{i+1}=z_{i}+\varepsilon \phi _{\xi _{i}\xi _{i+1}}(x_{i}(z_{i},\xi
,\varepsilon ),y_{i+1}(z_{i},\xi ,\varepsilon ),z_{i},\varepsilon ),
\label{eq:FullSlowDiscrete}
\end{equation}%
where $\phi \in C^{1}$, \ and the index $i$ counts intersections with the
Poincar\'{e} sections. The iterate $z_{i}\in \Sigma _{\xi _{i}}.$ Note that
the functions $x_{i}$ and $y_{i}$ depends on entire code $\xi .$

Lemma 2 in \cite{BG2007} implies that for $\varepsilon =0$ and any two codes 
$\xi ^{(1)}$ and $\xi ^{(2)}$ that satisfy $\xi _{i}^{(1)}=\xi _{i}^{(2)}$
for $\left\vert i\right\vert \leq n$ for any $n$ the following estimate
holds 
\begin{equation}
\max \left\{ \left\Vert x_{i}(z,\xi ^{(1)})-x_{i}(z,\xi ^{(2)})\right\Vert
,\left\Vert y_{i}(z,\xi ^{(1)})-y_{i}(z,\xi ^{(2)})\right\Vert \right\} \leq
2r\lambda ^{n-\left\vert i\right\vert },  \label{eq:DistCodes}
\end{equation}%
where the constants $r>0$ and $0<\lambda <1$ do not depend on the sequences $%
\xi ^{(1),(2)}.$ For $\varepsilon $ small Lemma 3 in \cite{BG2007} implies
the functions ($x_{i},y_{i}$) are defined for all small $\varepsilon $ and
all $z\in \overline{{{\mathcal{D}}}},$ they are uniformly bounded along with
their first derivatives with respect to $z$ and satisfy (\ref{eq:DistCodes}%
). Moreover, by the lemma, there is a constant $C_{0}>0$, independent of the
code $\xi ,$ such that 
\begin{equation}
\left\Vert x_{i}(z,\xi ,\varepsilon )-x_{i}(z,\xi ,0),y_{i}(z,\xi
,\varepsilon )-y_{i}(z,\xi ,0)\right\Vert <C_{0}\varepsilon ,
\label{eq:DistNH_Mflds}
\end{equation}%
for all $i\in {{\mathbb{Z}}}.$

\section{Proof of Theorem\label{sec:ProofMainThm} \protect\ref%
{Theorem:NewMain}}

In this section we construct a code $\xi ^{(\ast )}$ that generates a
trajectory $z^{(\ast )}$ of the full dynamics (\ref{eq:FullSlowDiscrete})
which satisfies Theorem \ref{Theorem:NewMain}. The proof relies on two
lemmas. The first one states that two trajectories of (\ref%
{eq:FullSlowDiscrete}) which have a slow iterate in common stay uniformly
close to each other as long as their respective codes coincide (Lemma \ref%
{Lemma:UniformBound}). The second one states that given any two points in
the slow phase space which are $O(\varepsilon )$ close we can find a
trajectory which goes from the neighbourhood of the first point to the
neighbourhood of the second point. This trajectory is essentially obtained
by updating the code of another trajectory. We then combine these two lemmas
to give an inductive proof of Theorem \ref{Theorem:NewMain}.

\subsection{Uniform closeness of trajectories}

As mentioned in the previous paragraph, finding a trajectory that satisfies
Theorem \ref{Theorem:NewMain} is an iterative process which involves
updating the code. The crucial point here is that when we update the code we
do not only alter the future of the trajectory but we switch to another
trajectory. Differently put, as we update the code the \emph{entire} slow
dynamics changes (not only the "future" iterates of the trajectory); in fact
if one considers the description of the slow component of the full dynamics
given in equation (\ref{eq:FullSlowDiscrete}) one sees that the first two
arguments of the function $\phi _{\xi _{i}\xi _{i+1}}$ depend on the \emph{%
entire} code $\xi $ (not just the current and next to current code
elements). The following lemma gives us a uniform estimate on how much the
"past" of the trajectory changes when the code is updated
and the initial condition is kept fixed, i.e. $z_{0}^{(a)}=z_{0}^{(b)}$. The
estimate appears in \cite{GT2007} (see the paragraph between equation (56)
and (57) therein) but since the proof was only sketched and the result is
crucial to our theory we state it together with a full proof.{}

\begin{lemma}
\label{Lemma:UniformBound} There exists $\varepsilon _{0}>0$ such that for
any $0<\varepsilon <\varepsilon _{0}$ the following holds: If $\xi ^{(a)}$
and $\xi ^{(b)}$ are two codes such that for some $N>0$%
\begin{equation*}
\xi _{i}^{(a)}=\xi _{i}^{(b)}\textrm{ \ for all }\left\vert i\right\vert \leq N
\end{equation*}%
and $z^{(a)}$ and $z^{(b)}$ are two slow trajectories of (\ref%
{eq:FullSlowDiscrete}) which correspond to the codes $\xi ^{(a)}$ and $\xi
^{(b)}$ respectively and $z_{0}^{(a)}=z_{0}^{(b)}$ then%
\begin{equation}
\left\Vert z_{i}^{(a)}-z_{i}^{(b)}\right\Vert \leq K\varepsilon \quad {{%
\textrm{for all}}}\quad 0\leq i\leq N\,  \label{Eq:UniformBound}
\end{equation}%
where 
\begin{equation}
K:=\frac{8\left\Vert \phi \right\Vert _{C^{1}}r\lambda }{1-\lambda }.
\label{Eq:UniformConst!}
\end{equation}
\end{lemma}

Note that the constant $K$ is independent of $N.$ We postpone the proof of
this Lemma to Section \ref{Sec:UniformBound}.

With $K$ given by (\ref{Eq:UniformConst!}) we let%
\begin{equation*}
A_{1}=\max_{c}\sup_{(u,v)\in \overline{\mathcal{D}}}3\varepsilon
(K+\left\Vert \phi _{cc}\right\Vert _{C^{1}}\frac{4r}{1-\lambda }%
+1+\left\vert T_{c}(u,v)\right\vert \left\vert X_{c}(u,v)\right\vert )
\end{equation*}%
and define recursively%
\begin{equation*}
A_{i}:=\max_{c}\sup_{(u,v)\in \overline{\mathcal{D}}}3\varepsilon
(A_{i-1}+\left\Vert \phi _{cc}\right\Vert _{C^{1}}\frac{4r}{1-\lambda }%
+1+\left\vert T_{c}(u,v)\right\vert \left\vert X_{c}(u,v)\right\vert ).
\end{equation*}%
Let%
\begin{equation}
A:=A_{2d}  \label{eq:A}
\end{equation}%
where $d$ is the number of slow degrees of freedom.

\begin{lemma}
\label{Lemma:Shadowing} Let $L>0$ be any fixed constant and $A$ given by (%
\ref{eq:A}). There exists $\varepsilon _{0}>0$ such that for any $%
0<\varepsilon <\varepsilon _{0}$ and any $z\in \overline{\mathcal{D}}$ and
any code $\xi ^{(a)}$ with a corresponding trajectory $z_{i}^{(a)}$
satisfying%
\begin{equation*}
\left\Vert z_{0}^{(a)}-z\right\Vert \leq \varepsilon (L+A)
\end{equation*}%
the following holds: There exists another code $\xi ^{(b)}$ with%
\begin{equation*}
\xi _{i}^{(b)}=\xi _{i}^{(a)}\textrm{ \ \ }i<0
\end{equation*}%
such that for any $p\leq 0$ there exists $N\in \mathbb{N}$ and a
corresponding trajectory $z^{(b)}$  which satisfies%
\begin{equation*}
z_{p}^{(a)}=z_{p}^{(b)},
\end{equation*}%
\begin{equation*}
\left\Vert z_{N}^{(b)}-z\right\Vert \leq \varepsilon A.
\end{equation*}%
Moreover there exists a uniformly bounded constant $C_{1}>0$ such that%
\begin{equation}
\left\Vert z_{i}^{(b)}-z\right\Vert \leq \varepsilon (C_{1}+A)\textrm{ \ for }%
0\leq i\leq N.  \label{eq:ShadowingBall}
\end{equation}
\end{lemma}

Next we combine these two lemmas to show that we can shadow any curve $%
\gamma :\mathbb{R}^{+}\rightarrow \overline{\mathcal{D}}$.

\subsection{\protect\bigskip Combining the results}

\begin{proof}
\textbf{(of Theorem \ref{Theorem:NewMain}) }Let $\gamma :\mathbb{R}%
^{+}\rightarrow \overline{\mathcal{D}}$ be any curve. Pick a constant $L>0$
arbitrarily. Take the smallest $\varepsilon _{0}$ of Lemma \ref%
{Lemma:UniformBound} and Lemma \ref{Lemma:Shadowing}. The following analysis
is valid for any $0<\varepsilon <\varepsilon _{0}$. Using $L$ we define a
sequence $t_{i}$, $i\in \mathbb{N}$, as follows%
\begin{equation}
\left\{ 
\begin{array}{l}
t_{0}=0 \\ 
\\ 
t_{i+1}=\min\limits_{t>t_{i}}\left\{ t:\left\Vert \gamma (t)-\gamma
(t_{i})\right\Vert =\varepsilon L\right\}%
\end{array}%
\right.  \label{eq:Curve_discrete}
\end{equation}%
and if for some $k$%
\begin{equation*}
\left\Vert \gamma (t)-\gamma (t_{k})\right\Vert \leq \varepsilon L\textrm{ \ \
for all }t>t_{k}
\end{equation*}%
then%
\begin{equation*}
t_{i+1}=t_{i}+1\textrm{ \ for all }i\geq k.
\end{equation*}%
The sequence $t_{i}$ divide the curve $\gamma $ into points $\gamma (t_{i})$
which are at most $\varepsilon L$ apart. Next we take any code $\xi ^{(0)}$
and a corresponding trajectory $z^{(0)}$ such that $z_{0}^{(0)}=\gamma
(t_{0})=\gamma (0).$ We also define $P(0)=0.$

\textbf{Inductive assumption: }There exists a code $\xi ^{(l)}$ and a
monotone sequence $P(l)$ such that for all $1\leq l\leq k$%
\begin{equation*}
\xi _{j}^{(l)}=\xi _{j}^{(l-1)}\textrm{ \ for }j<P(l-1).
\end{equation*}%
Moreover, there is a trajectory $z^{(l)}$ corresponding to $\xi ^{(l)}$ such
that%
\begin{equation*}
z_{0}^{(l)}=z_{0}^{(0)}
\end{equation*}%
and%
\begin{equation}
\left\Vert z_{P(l)}^{(l)}-\gamma (t_{l})\right\Vert \leq \varepsilon A.
\label{eq:IndHyp1}
\end{equation}%
Furthermore%
\begin{equation}
\left\Vert z_{j}^{l}-\gamma (t_{l})\right\Vert \leq \varepsilon (C_{1}+A)%
\textrm{ \ for \ }P(l-1)\leq j\leq P(l),  \label{eq:IndHyp2}
\end{equation}%
where $C_{1}$ is given by Lemma \ref{Lemma:Shadowing}.

\textbf{Inductive step.} We will use Lemma \ref{Lemma:Shadowing} to verify
the inductive assumption. Set $p=-P(k),$ $z=\gamma (t_{k+1})$ and $%
z_{i}^{(a)}=z_{P(k)+i}^{(k)}$ then using (\ref{eq:IndHyp1}) and (\ref%
{eq:Curve_discrete}) we get%
\begin{equation*}
\left\Vert z_{0}^{(a)}-z\right\Vert \leq \varepsilon (A+L).
\end{equation*}%
Then applying Lemma \ref{Lemma:Shadowing} implies that there exists a code $%
\xi _{P(k)+i}^{(k+1)}:=\xi _{i}^{(b)}$ which satisfies%
\begin{equation*}
\xi _{i}^{(k+1)}=\xi _{i}^{(k)}\textrm{ \ for }i<P(k)
\end{equation*}%
and generates a trajectory $z_{j}^{(k+1)}$ which satisfies%
\begin{equation*}
z_{0}^{(k+1)}=z_{0}^{(k)}=\gamma (0).
\end{equation*}%
Furthermore there exists $N_{k+1}$ such that%
\begin{equation*}
\left\Vert z_{P(k)+N_{k+1}}^{(k+1)}-z\right\Vert \leq \varepsilon A
\end{equation*}%
and%
\begin{equation*}
\left\Vert z_{j}^{k+1}-\gamma (t_{k+1})\right\Vert \leq \varepsilon (C+A)%
\textrm{ \ for }P(k)\leq j\leq P(k+1).
\end{equation*}%
The monotone sequence $P(k)$ is defined inductively by%
\begin{eqnarray*}
P(0) &=&0, \\
P(k) &=&P(k-1)+N_{k}.
\end{eqnarray*}%
This concludes the inductive step.

By induction there exists a unique code $\xi ^{(\ast )}$ which for all $k>0$
satisfies%
\begin{equation*}
\xi _{i}^{(\ast )}=\xi _{i}^{(k)}\textrm{ \ \ for }i\leq P(k).
\end{equation*}%
Denote by $z^{(\ast )}$ the trajectory of (\ref{eq:FullSlowDiscrete}) which
corresponds to the code $\xi ^{(\ast )}$ and satisfies $z_{0}^{(\ast
)}=\gamma (0).$ Next we want to apply Lemma \ref{Lemma:UniformBound} to show
that for all $k$ the trajectory $z^{(\ast )}$ lies close to the points $%
\gamma (t_{k}).$ Set $\xi ^{(a)}=\xi ^{(\ast )},$ $\xi ^{(b)}=\xi ^{(k)}$
and $N=P(k).$ Then Lemma \ref{Lemma:UniformBound} implies%
\begin{equation*}
\left\Vert z_{i}^{(\ast )}-z_{i}^{(k)}\right\Vert \leq \varepsilon K\textrm{ \
for }0\leq i\leq P(k).
\end{equation*}%
Combining this estimate with (\ref{eq:IndHyp2}) gives%
\begin{equation}
\forall k>0\textrm{ \ \ \ }\left\Vert z_{i}^{(\ast )}-\gamma
(t_{k})\right\Vert \leq \varepsilon (K+C_{1}+A)\textrm{ \ for }P(k-1)\leq
i\leq P(k).  \label{eq:z*-curve}
\end{equation}

To conclude the proof we have to define the time reparameterisation $%
\mathcal{T}(t)$ of Theorem \ref{Theorem:NewMain} and pass from the discrete
solution $z_{i}^{(\ast )}$ to the time continuous $z^{(\ast )}(t)$. Let us
begin with the time reparameterisation. For $t_{k}\leq t\leq t_{k+1}$ we
define%
\begin{equation*}
\mathcal{T}_{k}(t)=\tau _{P(k)}+\frac{\tau _{P(k+1)}-\tau _{P(k)}}{%
t_{k+1}-t_{k}}(t-t_{k}),
\end{equation*}%
where $\tau _{i}$ is the time $z^{(\ast )}(t)$ intersects the Poincar\'{e}
surface $\Sigma _{\xi _{i}^{(\ast )}},$ i.e. $z^{(\ast )}(\tau
_{i})=z_{i}^{(\ast )}.$ The complete $\mathcal{T}(t)$ is obtained by gluing
together all $\mathcal{T}_{k}(t).$

Next for $t_{k}\leq t\leq t_{k+1}$, and consequently $P(k)\leq i\leq P(k+1),$
we estimate the distance between $z^{(\ast )}(t)$ and $\gamma (t)$ as follows%
\begin{eqnarray*}
&&\left\Vert z^{(\ast )}(\mathcal{T}(t))-\gamma (t)\right\Vert \\
&\leq &\left\Vert z^{(\ast )}(\mathcal{T}(t))-z^{(\ast )}(\mathcal{T}%
(t_{k+1}))\right\Vert +\left\Vert z^{(\ast )}(\mathcal{T}%
(t_{k+1}))-z_{P(k+1)}^{(\ast )}\right\Vert \\
&&+\left\Vert z_{P(k+1)}^{(\ast )}-\gamma (t_{k+1})\right\Vert +\left\Vert
\gamma (t_{k+1})-\gamma (t)\right\Vert .
\end{eqnarray*}%
Consider the right hand side. The second term is $0$ by definition, the
third term is bounded by (\ref{eq:z*-curve}) and the fourth one by (\ref%
{eq:Curve_discrete}). The first term in the right hand side we estimate as%
\begin{equation*}
\left\Vert z^{(\ast )}(\mathcal{T}(t))-z^{(\ast )}(\mathcal{T}%
(t_{k+1}))\right\Vert \leq \sup_{(u,v)\in \overline{\mathcal{D}}}\dot{z}%
\sup_{k}(\tau _{P(k+1)}-\tau _{P(k)}).
\end{equation*}%
We note that the set of periodic orbits and transversal heteroclinics is
compact. The full trajectory, whose projection on the slow phase space is $%
z^{(\ast )},$ lies in a compact neighbourhood of this set, therefore $\dot{z}
$ given by (\ref{eq:FullDyn}) is uniformly bounded on this set: there exists
a constant $C_{2}>0,$ independent of $k,$ such that%
\begin{equation*}
\left\Vert z^{(\ast )}(\mathcal{T}(t))-z^{(\ast )}(\mathcal{T}%
(t_{k+1}))\right\Vert \leq \varepsilon C_{2}.
\end{equation*}%
By choosing $C:=C_{2}+L+K+C_{1}+A$ we have%
\begin{equation*}
\left\Vert z^{(\ast )}(\mathcal{T}(t))-\gamma (t)\right\Vert \leq
\varepsilon C\textrm{ \ for }t>0.
\end{equation*}
\end{proof}

\section{Technical Results}

\label{Sec:Technical}

In this section we include the proofs of Lemma \ref{Lemma:UniformBound} and
Lemma \ref{Lemma:Shadowing} that we used to prove Theorem \ref%
{Theorem:NewMain}

\subsection{Proof of Lemma \protect\ref{Lemma:UniformBound}\label%
{Sec:UniformBound}}

Let $\varepsilon _{0}$ be sufficiently small for the normal hyperbolicity
estimates to be valid and smaller than $\frac{1-\lambda }{2C\lambda },$
where 
\begin{equation}
C=\sup_{(u,v)\in \overline{\mathcal{D}}}\left\Vert \phi \right\Vert
_{C^{1}}\left( 1+\max \left( \left\Vert \frac{\partial x}{\partial z}%
\right\Vert ,\left\Vert \frac{\partial y}{\partial z}\right\Vert \right)
\right) .  \label{eq:C}
\end{equation}%
Consider the two codes $\xi ^{(a)}$ and $\xi ^{(b)}$. By assumption we have 
\begin{equation*}
\xi _{i}^{(a)}=\xi _{i}^{(b)}\quad \left\vert i\right\vert \leq N,
\end{equation*}%
and as a consequence of normal hyperbolicity that 
\begin{equation}
\left\Vert x_{i}(z,\xi ^{(a)})-x_{i}(z,\xi ^{(b)}),y_{i}(z,\xi
^{(a)})-y_{i}(z,\xi ^{(b)})\right\Vert \leq 2r\lambda ^{N-\left\vert
i\right\vert },  \label{Eq:ExpCloseness}
\end{equation}%
for all $\left\vert i\right\vert \leq N,$ see (\ref{eq:DistCodes})$.$ For
convenience we repeat equation (\ref{eq:FullSlowDiscrete}) which gives the
slow component of the full dynamics generated by the code $\xi $ 
\begin{equation*}
z_{i+1}=z_{i}+\varepsilon \phi _{\xi _{i}\xi _{i+1}}(x_{i}(z_{i},\varepsilon
;\xi ),y_{i+1}(z_{i},\varepsilon ;\xi ),z_{i},\varepsilon ).
\end{equation*}%
We note that the functions $x_{i}$ and $y_{i}$ depend on the full code $\xi $
(not only the current element in the code)$.$ Taking the difference of the $%
z $ components of the two trajectories coded by $\xi ^{(a)}$ and $\xi ^{(b)}$
and using that they have the same initial condition $z_{0}^{(a)}=z_{0}^{(b)}$
we get 
\begin{eqnarray*}
\left\Vert z_{k}^{(a)}-z_{k}^{(b)}\right\Vert &=&\varepsilon \left\Vert
\sum_{i=0}^{k-1}\phi _{\xi _{i}^{(a)}\xi
_{i+1}^{(a)}}(x_{i}(z_{i}^{(a)},\varepsilon ;\xi
^{(a)}),y_{i+1}^{{}}(z_{i}^{(a)},\varepsilon ;\xi
^{(a)}),z_{i}^{(a)},\varepsilon )-\right. \\
&&\left. \phi _{\xi _{i}^{(b)}\xi
_{i+1}^{(b)}}(x_{i}(z_{i}^{(b)},\varepsilon ;\xi
^{(b)}),y_{i+1}^{{}}(z_{i}^{(b)},\varepsilon ;\xi
^{(b)}),z_{i}^{(b)},\varepsilon )\right\Vert .
\end{eqnarray*}%
Now, $\xi _{i}^{(a)}=\xi _{i}^{(b)}$ for all $\left\vert i\right\vert <N,$
therefore $\phi _{\xi _{i}^{(a)}\xi _{i+1}^{(a)}}$ and $\phi _{\xi
_{i}^{(b)}\xi _{i+1}^{(b)}}$ are the same functions in the interval we are
studying. Using the mean value inequality and suppressing the dependence of $%
\varepsilon $ in the notation we get 
\begin{eqnarray*}
\left\Vert z_{k}^{(a)}-z_{k}^{(b)}\right\Vert &\leq &\varepsilon \left\Vert
\phi \right\Vert _{C^{1}}\left\Vert \sum_{i=0}^{k-1}\left(
x_{i}(z_{i}^{(a)},\xi ^{(a)})-x_{i}(z_{i}^{(b)},\xi ^{(b)}),\right. \right.
\\
&&\left. \left. y_{i+1}(z_{i}^{(a)},\xi ^{(a)})-y_{i+1}(z_{i}^{(b)},\xi
^{(b)})\right) \right\Vert + \\
&&\varepsilon \left\Vert \phi \right\Vert _{C^{1}}\left\Vert
\sum_{i=0}^{k-1}(z_{i}^{(a)}-z_{i}^{(b)})\right\Vert .
\end{eqnarray*}%
By adding and subtracting terms and using the triangle inequality we rewrite
this expression to a form where inequality $\left( \ref{Eq:ExpCloseness}%
\right) $ can be used 
\begin{eqnarray*}
\left\Vert z_{k}^{(a)}-z_{k}^{(b)}\right\Vert &\leq &\varepsilon \left\Vert
\phi \right\Vert _{C^{1}}\sum_{i=0}^{k-1}\left\Vert \left(
x_{i}(z_{i}^{(a)},\xi ^{(a)})-x_{i}(z_{i}^{(b)},\xi ^{(b)}),\right. \right.
\\
&&\left. \left. y_{i+1}(z_{i}^{(a)},\xi ^{(a)})-y_{i+1}(z_{i}^{(a)},\xi
^{(b)})\right) \right\Vert + \\
&&\varepsilon \left\Vert \phi \right\Vert _{C^{1}}\sum_{i=0}^{k-1}\left\Vert
\left( x_{i}(z_{i}^{(a)},\xi ^{(b)})-x_{i}(z_{i}^{(b)},\xi ^{(b)}),\right.
\right. \\
&&\left. \left. y_{i+1}(z_{i}^{(a)},\xi ^{(b)})-y_{i+1}(z_{i}^{(b)},\xi
^{(b)})\right) \right\Vert + \\
&&\varepsilon \left\Vert \phi \right\Vert _{C^{1}}\sum_{i=0}^{k-1}\left\Vert
z_{i}^{(a)}-z_{i}^{(b)}\right\Vert .
\end{eqnarray*}%
Using inequality $\left( \ref{Eq:ExpCloseness}\right) $ to estimate the
first term and the mean value inequality to estimate the second term gives
us 
\begin{eqnarray*}
\left\Vert z_{k}^{(a)}-z_{k}^{(b)}\right\Vert &\leq &\varepsilon \left\Vert
\phi \right\Vert _{C^{1}}\sum_{i=0}^{k-1}2r\lambda ^{N-i}+\varepsilon
\left\Vert \phi \right\Vert _{C^{1}}\sum_{i=0}^{k-1}\left\Vert
z_{i}^{(a)}-z_{i}^{(b)}\right\Vert \\
&&+\varepsilon \left\Vert \phi \right\Vert _{C^{1}}\max \left( \left\Vert 
\frac{\partial x}{\partial z}\right\Vert ,\left\Vert \frac{\partial y}{%
\partial z}\right\Vert \right) \sum_{i=0}^{k-1}\left\Vert
z_{i}^{(a)}-z_{i}^{(b)}\right\Vert ,
\end{eqnarray*}%
which we rewrite as 
\begin{equation}
\left\Vert z_{k}^{(a)}-z_{k}^{(b)}\right\Vert \leq 2\varepsilon \left\Vert
\phi \right\Vert _{C^{1}}r\frac{\lambda ^{N-k+1}}{1-\lambda }+\varepsilon
C\sum_{i=0}^{k-1}\left\Vert z_{i}^{(a)}-z_{i}^{(b)}\right\Vert ,
\label{Eq:z_diff_Gronwall}
\end{equation}%
where $C$ is given by (\ref{eq:C}). Now, in order to prove Lemma \ref%
{Lemma:UniformBound} we will show that 
\begin{equation}
\left\Vert z_{k}^{(a)}-z_{k}^{(b)}\right\Vert \leq \varepsilon K\lambda
^{N-k},  \label{Eq:UniformInduction}
\end{equation}%
for all $0\leq k<N$ where $K$ is given by $\left( \ref{Eq:UniformBound}%
\right) .$ The Gronwall type of estimate follows from inequality $\left( \ref%
{Eq:z_diff_Gronwall}\right) $ by finite induction. Since $%
z_{0}^{(a)}=z_{0}^{(b)}$ the statement is valid for $k=0.$ Let us assume
that $\left( \ref{Eq:UniformInduction}\right) $ is true for all $k\leq m.$
Consider the case $k=m+1.$ By $\left( \ref{Eq:z_diff_Gronwall}\right) $ we
get 
\begin{equation*}
\left\Vert z_{m+1}^{(a)}-z_{m+1}^{(b)}\right\Vert \leq 2\varepsilon
\left\Vert \phi \right\Vert _{C^{1}}r\frac{\lambda ^{N-m}}{1-\lambda }%
+\varepsilon C\sum_{i=0}^{m}\left\Vert z_{i}^{(a)}-z_{i}^{(b)}\right\Vert .
\end{equation*}%
The induction assumption implies that 
\begin{eqnarray*}
\left\Vert z_{m+1}^{(a)}-z_{m+1}^{(b)}\right\Vert &\leq &2\varepsilon
\left\Vert \phi \right\Vert _{C^{1}}r\frac{\lambda ^{N-m}}{1-\lambda }%
+\varepsilon C\sum_{i=0}^{m}\varepsilon K\lambda ^{N-i} \\
&\leq &2\varepsilon \left\Vert \phi \right\Vert _{C^{1}}r\frac{\lambda ^{N-m}%
}{1-\lambda }+\varepsilon ^{2}CK\frac{\lambda ^{N-m}}{1-\lambda }.
\end{eqnarray*}%
For%
\begin{equation*}
K:=\frac{8\left\Vert \phi \right\Vert _{C^{1}}r\lambda }{1-\lambda }
\end{equation*}%
and $0<$ $\varepsilon <\frac{1-\lambda }{2C\lambda }$ we have 
\begin{equation*}
K>\frac{2\left\Vert \phi \right\Vert _{C^{1}}r\lambda }{1-\lambda
-\varepsilon C\lambda }.
\end{equation*}%
Thus it follows that 
\begin{equation}
\left\Vert z_{m+1}^{(a)}-z_{m+1}^{(b)}\right\Vert \leq \varepsilon K\lambda
^{N-(m+1)}.  \label{Eq:InductionComplete}
\end{equation}%
Since $0<\lambda <1$ and $0<m+1\leq N$ the lemma follows from $\left( \ref%
{Eq:InductionComplete}\right) .$ We note that $K$ is independent of $N$ and $%
\varepsilon .\hfill \square $

\subsection{Proof of Lemma \protect\ref{Lemma:Shadowing}}

Fix $L>0$ arbitrarily and let $A$ be given by (\ref{eq:A}). Then we choose $%
\varepsilon _{0}$ as the smallest of that required by Lemma \ref%
{Lemma:UniformBound} and%
\begin{equation*}
\frac{1}{C\max_{c}\left( \sup_{(u,v)\in \overline{\mathcal{D}}}\frac{L+A}{%
T_{c}(u,v)\left\vert \nabla J_{c}\right\vert }\right) },
\end{equation*}%
where $C$ is given by (\ref{eq:C}). Then for any $\varepsilon <\varepsilon
_{0}$ let $\xi ^{(a)}$ be a code generating a trajectory of (\ref%
{eq:FullSlowDiscrete}) such that%
\begin{equation*}
\left\Vert z_{0}^{(a)}-z\right\Vert \leq \varepsilon (L+A).
\end{equation*}%
To show that there exists a code $\xi ^{(b)}$ that satisfies the lemma we
introduce the notion of guiding path. Denoting by $T_{z_{0}^{(a)}}(\mathbb{R}%
^{2d})$ the tangent space of $\mathbb{R}^{2d}$ at the point $z_{0}^{(a)}$,
consider the vector $\overrightarrow{v}\in T_{z_{0}^{(a)}}(\mathbb{R}^{2d})$
pointing towards $z$ with length $\varepsilon (L+A).$ Condition $\left[ 
\mathbf{A3}\right] $ implies that there exists (a possibly non-unique)
injective map $\sigma :\left\{ 1,\ldots ,2d\right\} \rightarrow \left\{
c_{1},\ldots ,c_{n}\right\} $ (where $c_{i}$ is the index of the periodic
orbits), with $\sigma =\sigma (\overrightarrow{v})$ and non-negative
coefficients $a_{1},\ldots ,a_{2d}$ such that%
\begin{equation}
\overrightarrow{v}=\sum_{i=1}^{2d}\varepsilon a_{i}\overrightarrow{X}%
_{\sigma (i)}(z_{0}^{(a)}).  \label{eq:linearcomb}
\end{equation}%
We refer to the coefficients $a_{i}$ as \emph{guiding times}. We denote%
\footnote{%
At this point we abuse the notation by identifying vector fields and points
in $\mathbb{R}^{2d}$.} the \emph{guiding path} between $z_{0}^{(a)}$ and $z$
as%
\begin{equation}
G(t,z_{0}^{(a)},z)=z_{0}^{(a)}+\sum_{j=1}^{i}\varepsilon a_{j}%
\overrightarrow{X}_{\sigma (j)}(z_{0}^{(a)})+\varepsilon
(t-\sum_{j=1}^{i}a_{j})\overrightarrow{X}_{\sigma (i)}(z_{0}^{(a)})
\label{eq:GuidingPath}
\end{equation}%
for $\varepsilon \sum_{j=1}^{i}a_{j}<t<\varepsilon \sum_{j=1}^{i+1}a_{j}.$
Let us denote by $W_{\overrightarrow{v}}:=\left\{ \overrightarrow{X}_{\sigma
(1)}(z_{0}^{(a)}),\ldots ,\overrightarrow{X}_{\sigma
(2d)}(z_{0}^{(a)})\right\} $ the $2d\times 2d$ matrix whose columns are the $%
2d$ Hamiltonian vectors needed to represent $\overrightarrow{v}$ as a linear
combination of the type (\ref{eq:linearcomb}). Now, condition $\left[ 
\mathbf{A3}\right] $ implies that $\left\{ \overrightarrow{X}_{\sigma
(1)}(z_{0}^{(a)}),\ldots ,\overrightarrow{X}_{\sigma
(2d)}(z_{0}^{(a)})\right\} $ is a basis in $\mathbb{R}^{2d},$ thus $W_{%
\overrightarrow{v}}$ is invertible and from (\ref{eq:linearcomb}), denoting $%
\overrightarrow{a}=(a_{1},\ldots ,a_{2d}),$ we get%
\begin{equation}
\left\Vert \overrightarrow{a}\right\Vert \leq \left\Vert W_{\overrightarrow{v%
}}^{-1}\right\Vert \frac{\left\Vert \overrightarrow{v}\right\Vert }{%
\varepsilon }\leq D  \label{eq:D}
\end{equation}%
where $D=\sup_{z\in \overline{\mathcal{D}}}\sup_{\overrightarrow{\zeta }\in
T_{z_{0}^{(a)}}(\mathbb{R}^{2d})}\left\Vert W_{\overrightarrow{\zeta }%
}^{-1}\right\Vert (L+A)<\infty $ depends only on $L+A.$ Hence the guiding
times $a_{i}$ are uniformly bounded.

Without loss of generality
\footnote{The proof for $d>1$ is analogous to the case $d=1.$ The only difference is
that the guiding path will consist of more segments. Consequently the
constant $A$ depends on $d.$} 
we will prove the lemma for $d=1,$ and we
write $v=\varepsilon a_{1}X_{1}+\varepsilon a_{2}X_{2}.$ Then the guiding
path consists of two segments. We will show that these two segments can be
shadowed one at a time, hence the code $\xi ^{(b)}$ will be obtained by
updating the code $\xi ^{(a)}$ twice. The rules for updating the code are 
\begin{equation}
\left\{ \begin{array}{ll}
{\tilde \xi}_{k}^{(b)}=\xi_{k}^{(a)} &  \textrm{for}\: k\leq 0, \\
{\tilde \xi}_{k}^{(b)}=c_{1} & \textrm{for all}\: k>0,
\label{Eq:Preliminary}
\end{array} \right .
\end{equation}
and
\begin{equation}
\left\{ \begin{array}{ll}
\xi _{k}^{(b)}=\tilde{\xi }_{k}^{(b)}\quad {{\textrm{for}}}\quad k\leq
\left\lceil \frac{a_{1}}{T_{c_{1}}}\right\rceil , \\ 
\xi _{k}^{(b)}=c_{2}\quad {{\textrm{for all}}}\quad k>\left\lceil \frac{a_{1}}{%
T_{c_{1}}}\right\rceil ,%
\end{array}
\right .
\label{Eq:Updated}
\end{equation}%
where $\left\lceil \cdot \right\rceil $ denotes rounding up to the next
integer and $T_{c_{1}}$ denotes the period of the periodic orbit $L_{c_{1}}$
(also note that by (\ref{eq:D}) $a_{i}$ are of order $O(\varepsilon )$)$.$
We define $\widetilde{N}=\left\lceil \frac{a_{1}}{T_{c_{1}}}\right\rceil $
and $N^{b}=\left\lceil \frac{a_{2}}{T_{c_{2}}}\right\rceil $ and let $N=%
\widetilde{N}+N^{b}.$

We begin by updating the code according to (\ref{Eq:Preliminary}) which
gives us the code $\widetilde{\xi }^{(b)}.$ For any $p\leq 0$ we pick the
trajectory of (\ref{eq:FullSlowDiscrete}) corresponding to $\widetilde{\xi }%
^{(b)}$ such that $\widetilde{z}_{p}^{(b)}=z_{p}^{(a)}.$ Let us also
consider two auxiliary sequences, $\hat{z}_{k}^{(b)}$ which corresponds to
the code $\widetilde{\xi }^{(b)}$ and $z_{k}^{c_{1}}$ which corresponds to
the constant code $\left( c_{1}\right) ^{\infty }$. Both sequences satisfy
the initial condition 
\begin{equation*}
\hat{z}_{0}^{(b)}=z_{0}^{c_{1}}=z_{0}^{(a)}\,.
\end{equation*}%
By the definition (\ref{Eq:Preliminary}) we have $\widetilde{\xi }%
_{k}^{(b)}=c_{1}$ for all $k\geq 0$. Then

\begin{lemma}
\label{Lemma:SameCodes}(Lemma 5, \cite{BG2007}) For any $K_{0}>0,t_{0}>0,$
there is $\varepsilon _{0}>0$ such that for any $\left\vert \varepsilon
\right\vert <\varepsilon _{0}$ and any two codes $\xi ^{1}$ and $\xi ^{2}$
such that for some index $j$%
\begin{equation*}
\xi _{j+i}^{1}=\xi _{j+i}^{2}=c\quad 0\leq i\leq N_{0}(\varepsilon )\equiv
\left\lfloor \frac{t_{0}}{\varepsilon }\right\rfloor
\end{equation*}%
the inequality $\left\Vert z_{j}^{1}-z_{j}^{2}\right\Vert \leq \varepsilon
K_{0}$ implies%
\begin{equation*}
\left\Vert z_{j+N}^{1}-z_{j+N}^{2}\right\Vert \leq \varepsilon
C_{1}e^{\varepsilon NC_{2}}\quad 0\leq N\leq N_{0}(\varepsilon ),
\end{equation*}%
where%
\begin{equation*}
C_{1}=\left\Vert \frac{\partial \phi _{cc}}{\partial (x,y)}\right\Vert \frac{%
4r}{1-\lambda }+K_{0}
\end{equation*}%
and%
\begin{equation*}
C_{2}=\left\Vert \frac{\partial \phi _{cc}}{\partial z}\right\Vert
+\left\Vert \frac{\partial \phi _{cc}}{\partial (x,y)}\right\Vert \max
\left\{ \left\Vert \frac{\partial x_{c}}{\partial z}\right\Vert ,\left\Vert 
\frac{\partial y_{c}}{\partial z}\right\Vert \right\} .
\end{equation*}%
In the compact set $\overline{\mathcal{D}}$ we have%
\begin{equation*}
\left\Vert z_{j+N}^{1}-z_{j+N}^{2}\right\Vert \leq 3\varepsilon C_{1}\quad
0\leq N\leq \frac{1}{\varepsilon C_{2}}.
\end{equation*}
\end{lemma}

implies that 
\begin{equation}
\left\Vert \hat{z}_{k}^{(b)}-{z}_{k}^{c_{1}}\right\Vert \leq 3\varepsilon
\left( \left\Vert \phi _{c_{1}c_{1}}\right\Vert _{C^{1}}\frac{4r}{1-\lambda }%
\right)  \label{Eq:star}
\end{equation}%
for all $0\leq k\leq N.$ By Lemma \ref{Lemma:UniformBound} we have the
following bound 
\begin{equation}
\Vert \widetilde{z}_{0}^{(b)}-\hat{z}_{0}^{(b)}\Vert =\Vert \widetilde{z}%
_{0}^{(b)}-z_{0}^{(a)}\Vert \leq \varepsilon K,  \label{Eq:Jump0}
\end{equation}%
since $\widetilde{z}_{k}^{(b)}$ and $z_{k}^{(a)}$ share the same code for $%
k\leq 0$ and $\widetilde{z}_{p}^{(b)}=z_{p}^{(a)}.$ Now, since $\widetilde{z}%
^{(b)}$ and $\hat{z}^{(b)}$ are generated from the identically same code,
Lemma~\ref{Lemma:SameCodes} implies that 
\begin{equation}
\left\Vert \widetilde{z}_{k}^{(b)}-\hat{z}_{k}^{(b)}\right\Vert \leq
3\varepsilon \left( \left\Vert \phi _{c_{1}c_{1}}\right\Vert _{C^{1}}\frac{4r%
}{1-\lambda }+K\right)  \label{eq:starstar}
\end{equation}%
for all $0\leq k\leq N$. Combining the estimates (\ref{Eq:star}) and (\ref%
{eq:starstar}) we obtain 
\begin{eqnarray}
\left\Vert \widetilde{z}_{k}^{(b)}-{z}_{k}^{c_{1}}\right\Vert &\leq
&\left\Vert \widetilde{z}_{k}^{(b)}-\hat{z}_{k}^{(b)}\right\Vert +\left\Vert 
\hat{z}_{k}^{(b)}-{z}_{k}^{c_{1}}\right\Vert   \\
&\leq &3\varepsilon \left( 2\left\Vert \phi _{c_{1}c_{1}}\right\Vert _{C^{1}}%
\frac{4r}{1-\lambda }+K\right) \,.  \label{Eq:Jump}
\end{eqnarray}%
for all $0\leq k\leq \widetilde{N}.$ Now it remains to prove that $z^{c_{1}}$
stays close to the guiding path $G.$ By definition we have%
\begin{equation*}
G(t,z_{0}^{(a)},z)=z_{0}^{(a)}-t\Omega ^{-1}\nabla J_{c_{1}}(z_{0}^{(a)})%
\textrm{ \ for }0\leq t\leq \varepsilon a_{1},
\end{equation*}%
where $\Omega ^{-1}$ is the inverse of the symplectic matrix. The map (\ref%
{eq:FullSlowDiscrete}) takes a simple form for the constant code $%
(c_{1})^{\infty }$, see \cite{BG2007} and the proof of Lemma 1 therein, 
\begin{equation*}
z_{k}^{c_{1}}=z_{0}^{(a)}-\varepsilon kT_{c_{1}}\Omega ^{-1}\nabla
J_{c_{1}}(z_{0}^{(a)})+O(\varepsilon ^{2}).
\end{equation*}%
Therefore%
\begin{eqnarray}
\left\Vert z_{\widetilde{N}}^{c_{1}}-G(\varepsilon
a_{1},z_{0}^{(a)},z)\right\Vert &\leq &O(\varepsilon ^{2})+\varepsilon
\max_{c}\sup_{(u,v)\in \overline{\mathcal{D}}}\left\vert
T_{c}(u,v)\right\vert \left\vert X_{c}(u,v)\right\vert   \\
&\leq &\varepsilon (1+\max_{c}\sup_{(u,v)\in \overline{\mathcal{D}}%
}\left\vert T_{c}(u,v)\right\vert \left\vert X_{c}(u,v)\right\vert ),
\label{eq:ShadowGP}
\end{eqnarray}%
where the $O(\varepsilon ^{2})$ term is the difference of the two maps and
the second term compensates for the fact that $a_{1}$ in general is not a
multiple of $T_{c_{1}},$ this round off error is bounded by the strength of
the vector field times the largest period. Collecting the estimates we obtain%
\begin{eqnarray}
\left\Vert \widetilde{z}_{\widetilde{N}}^{(b)}-G(\varepsilon
a_{1},z_{0}^{(a)},z)\right\Vert &\leq &3\varepsilon \left( 2\left\Vert \phi
_{c_{1}c_{1}}\right\Vert _{C^{1}}\frac{4r}{1-\lambda }+K\right) +\varepsilon
\label{eq:ShadowG_part1} \\
&&+\varepsilon \max_{c}\sup_{(u,v)\in \overline{\mathcal{D}}}\left\vert
T_{c}(u,v)\right\vert \left\vert X_{c}(u,v)\right\vert ,
\end{eqnarray}%
which is the accuracy with which we have shadowed the first segment of the
guiding path.

Let us continue by shadowing the second segment of the guiding path. We
begin by updating the code to $\xi ^{(b)},$ using the rule given by equation
(\ref{Eq:Updated}). We pick the trajectory $z^{(b)}$ which corresponds to $%
\xi ^{(b)}$ and satisfies $z_{p}^{(b)}=z_{p}^{(a)}$. We repeat the arguments
above to shadow the second segment of the guiding path. Indeed, we consider
two auxiliary sequences $\check{z}_{k}^{(b)}$ which corresponds to the code $%
\xi ^{(b)}$ and $z_{k}^{c_{2}}$ which corresponds to the constant code $%
\left( c_{2}\right) ^{\infty }$. Both sequences satisfy the initial
condition 
\begin{equation*}
\check{z}_{\widetilde{N}}^{(b)}=z_{\widetilde{N}}^{c_{2}}=G(\varepsilon
a_{1},z_{0}^{(a)},z)\,.
\end{equation*}%
By the definition (\ref{Eq:Preliminary}) we have $\xi _{k}^{(b)}=c_{2}$ for
all $k>\widetilde{N}$. Then by Lemma \ref{Lemma:SameCodes} 
\begin{equation}
\left\Vert \check{z}_{\widetilde{N}+k}^{(b)}-{z}_{\widetilde{N}%
+k}^{c_{2}}\right\Vert \leq 3\varepsilon \left( \left\Vert \phi
_{c_{2}c_{2}}\right\Vert _{C^{1}}\frac{4r}{1-\lambda }\right)
\label{Eq:star'}
\end{equation}%
for all $0\leq k\leq N^{b}.$ Lemma \ref{Lemma:UniformBound} together with (%
\ref{eq:ShadowG_part1}) yields the following bound 
\begin{eqnarray*}
\Vert z_{\widetilde{N}}^{(b)}-\check{z}_{\widetilde{N}}^{(b)}\Vert &=&\Vert
z_{\widetilde{N}}^{(b)}-G(\varepsilon a_{1},z_{0}^{(a)},z))\Vert \\
&\leq &\max_{c}\sup_{(u,v)\in \overline{\mathcal{D}}}(6\varepsilon
\left\Vert \phi _{cc}\right\Vert _{C^{1}}\frac{4r}{1-\lambda }+4\varepsilon K
\\
&&+\varepsilon +\left. \varepsilon \left\vert T_{c}(u,v)\right\vert
\left\vert X_{c}(u,v)\right\vert \right) .
\end{eqnarray*}%
Then Lemma~\ref{Lemma:SameCodes} implies that 
\begin{eqnarray}
&&\left\Vert z_{\widetilde{N}+k}^{(b)}-\check{z}_{\widetilde{N}%
+k}^{(b)}\right\Vert  \label{Eq:starstar'} \\
&\leq &3\varepsilon \max_{c}\sup_{(u,v)\in \overline{\mathcal{D}}}\left(
7\left\Vert \phi _{cc}\right\Vert _{C^{1}}\frac{4r}{1-\lambda }%
+4K+1+\left\vert T_{c}(u,v)\right\vert \left\vert X_{c}(u,v)\right\vert
\right)  
\end{eqnarray}%
for all $0\leq k\leq N^{b}$. Combining the estimates (\ref{Eq:star'}) and (%
\ref{Eq:starstar'}) gives 
\begin{eqnarray}
&&\left\Vert z_{\widetilde{N}+k}^{(b)}-{z}_{\widetilde{N}+k}^{c_{2}}\right%
\Vert \leq \left\Vert z_{\widetilde{N}+k}^{(b)}-\check{z}_{\widetilde{N}%
+k}^{(b)}\right\Vert +\left\Vert \check{z}_{\widetilde{N}+k}^{(b)}-{z}_{%
\widetilde{N}+k}^{c_{2}}\right\Vert  \label{eq:dagger} \\
&\leq &3\varepsilon \max_{c}\sup_{(u,v)\in \overline{\mathcal{D}}}\left(
8\left\Vert \phi _{cc}\right\Vert _{C^{1}}\frac{4r}{1-\lambda }%
+4K+1+\left\vert T_{c}(u,v)\right\vert \left\vert X_{c}(u,v)\right\vert
\right) \,.  
\end{eqnarray}%
for all $0\leq k\leq N^{b}.$ Using that $z_{\widetilde{N}%
}^{c_{2}}=G(t^{1},z_{0}^{(a)},z)$ the guiding path is shadowed by $z_{%
\widetilde{N}+k}^{c_{2}}$ exactly analogously to the first segment. The
result is%
\begin{equation}
\left\Vert z_{\widetilde{N}+k}^{c_{2}}-G(\varepsilon
a_{2},z_{0}^{(a)},z)\right\Vert \leq \varepsilon +\varepsilon
\max_{c}\sup_{(u,v)\in \overline{\mathcal{D}}}\left\vert
T_{c}(u,v)\right\vert \left\vert X_{c}(u,v)\right\vert  \label{eq:dagger'}
\end{equation}%
Using that $z_{\widetilde{N}+N^{b}}^{(b)}=z_{N}^{(b)}$ by definition and
combining (\ref{eq:dagger}) and (\ref{eq:dagger'}) gives%
\begin{eqnarray}
&&\left\Vert z_{N}^{(b)}-G(\varepsilon a_{2},z_{0}^{(a)},z)\right\Vert =
\label{eq:ShadowG_part2} \\
&\leq &4\varepsilon \max_{c}\sup_{(u,v)\in \overline{\mathcal{D}}}\left(
6\left\Vert \phi _{cc}\right\Vert _{C^{1}}\frac{4r}{1-\lambda }%
+3K+1+\left\vert T_{c}(u,v)\right\vert \left\vert X_{c}(u,v)\right\vert
\right)  
\end{eqnarray}

Thus, using the definition (\ref{eq:A}) of $A$%
\begin{equation*}
A=4\max_{c}\sup_{(u,v)\in \overline{\mathcal{D}}}\left( 6\left\Vert \phi
_{cc}\right\Vert _{C^{1}}\frac{4r}{1-\lambda }+3K+1+\left\vert
T_{c}(u,v)\right\vert \left\vert X_{c}(u,v)\right\vert \right) ,
\end{equation*}%
and the definition of guiding path gives%
\begin{equation*}
\left\Vert z_{N}^{(b)}-z\right\Vert \leq \varepsilon A.
\end{equation*}%
The constant $A$ is uniformly bounded on the compact space $\overline{%
\mathcal{D}}$. Lastly, from the fact that the guiding times $a_{i}$ are
uniformly bounded it follows that the length of the guiding path is bounded
by $\varepsilon C_{1}$ where $C_{1}$ is a uniform constant. This and the
estimates above implies that 
\begin{equation*}
\left\Vert z_{k}^{(b)}-z\right\Vert \leq \varepsilon (C_{1}+A)\textrm{ \ for }%
0\leq k\leq N.
\end{equation*}

\end{document}